\documentclass[12pt]{article}
\usepackage{amssymb}
\usepackage{epsfig}
\usepackage{graphicx}

\begin{document}

\newtheorem{mthm}{Theorem}
\newtheorem{mcor}{Corollary}
\newtheorem{mpro}{Proposition}
\newtheorem{mfig}{figure}
\newtheorem{mlem}{Lemma}
\newtheorem{mdef}{Definition}
\newtheorem{mrem}{Remark}
\newtheorem{mpic}{Picture}
\newtheorem{rem}{Remark}[section]
\newcommand{\ra}{{\mbox{$\rightarrow$}}}
\newtheorem{thm}{Theorem}[section]
\newtheorem{pro}{Proposition}[section]
\newtheorem{lem}{Lemma}[section]
\newtheorem{defi}{Definition}[section]
\newtheorem{cor}{Corollary}[section]

\title{\bf Asymptotic behavior of positive solutions to a degenerate elliptic equation in the upper half space with a nonlinear boundary condition}

\author{ Zhuoran Du }
\footnotetext{College of Mathematics and Econometrics, Hunan University, Changsha 410082,
                 PRC.
 E-mail: {duzr@hnu.edu.cn}}

\date{}
\maketitle

\begin{abstract}
We consider  positive  solutions of the problem
\begin{equation}
\left\{\begin{array}{l}-\mbox{div}(x_{n}^{a}\nabla u)=0\;\;\;\;\;\mbox{in}\;\;\mathbb{R}_{+}^{n},\\
\frac{\partial u}{\partial \nu^a}=u^{q} \;\;\;\;\;\;\;\;
\;\;\;\;\;\;\;\;\;\mbox{on}\;\;\partial \mathbb{R}_{+}^{n},\\
\end{array}
\right.
\end{equation}
where  $a\in (-1,0)\cup(0,1)$,  $q>1$ and
$\frac{\partial u}{\partial \nu^a}:=-\lim_{x_{n}\rightarrow
0^+}x_{n}^{a}\frac{\partial u}{\partial x_{n}}$.
 We  obtain  some qualitative properties of positive axially symmetric solutions   in
 $n\geq3$ for the case $a\in
(-1,0)$ under the condition $q\geq\frac{n-a}{n+a-2}$.
In particular, we establish the asymptotic expansion of positive axially symmetric solutions.

\end{abstract}



\section{Introduction}

In this paper, we consider qualitative properties of  positive  solutions of the problem
\begin{equation}\label{e2}
\left\{\begin{array}{l}-\mbox{div}(x_{n}^{a}\nabla u)=0\;\;\;\;\;\mbox{in}\;\;\mathbb{R}_{+}^{n}=\{x\in \mathbb{R}^{n}: x_{n}>0\},\\
\frac{\partial u}{\partial \nu^a}=u^{q} \;\;\;\;\;\;\;\;
\;\;\;\;\;\;\;\;\;\mbox{on}\;\;\partial \mathbb{R}_{+}^{n}=\{x\in \mathbb{R}^{n}: x_{n}=0\},\\
\end{array}
\right.
\end{equation}
where  $a\in (-1,0)\cup(0,1)$,  $q>1$,  $x=(x',x_{n})\in
\mathbb{R}^{n-1}\times [0,+\infty)$,   and $\frac{\partial u}{\partial \nu^a}:=-\lim_{x_{n}\rightarrow
0^+}x_{n}^{a}\frac{\partial u}{\partial x_{n}}$.

When the parameter
$a=0$, there had been many results about existence and  qualitative properties of
positive  solutions to the problem (\ref{e2}).  For critical case $q=\frac{n}{n-2}$, positive
solutions are known to exist in dimensions $n>2$. Moreover, any
positive solution is of the form (see \cite{MMM1} and \cite{LZ})
$$ u(x)
=\alpha|x-x^0|^{-(n-2)}, ~~\alpha>0, ~~x^0=(x^0_1, \ldots, x^0_n)\in \mathbb{R}^n,
~~x^0_n=-\frac{1}{n-2}\alpha^{\frac{2}{n-2}}.
$$
Positive solutions do not exist for subcritical case $q<\frac{n}{n-2}$
(see \cite{HU}). Subsequently, Chipot-Chlebik-Fila-Shafrir \cite{MMM2} obtained the existence of positive solutions of
(\ref{e2}) for the supercritical case $q>
\frac{n}{n-2}$. Recently,
under the conditions  $n>2$ and $q>\frac{n}{n-2}$, Harada \cite{H} established the
asymptotic expansion and the intersection property of positive
solutions to the problem (\ref{e2}).
Qui-Reichel  \cite{Qui} proved the existence of a unique singular positive axial symmetric  solution 
for   $q\geq\frac{n-1}{n-2}$,
in $n\geq3$. 

It is known that the problem (\ref{e2}) with $a=0$ is related to the square root of Laplacian  equation
\begin{equation}\label{e3}
(-\Delta)^{\frac{1}{2}}v=v^{q},\;\;\;\;\mbox{in}\;\mathbb{R}^{n-1}.
\end{equation}
Indeed, from the well-known Caffarelli-Silvestre extension in \cite{CS2} we know that if $u$ is a positive  solution of (\ref{e2}), then a
positive constant multiple of $u(x',0)=:v(x')$ satisfies (\ref{e3}).
Therefore the
asymptotic expansion results in \cite{H} generalize the asymptotic expansion results of the corresponding semilinear elliptic equation
with standard Laplacian operator
$$
-\Delta v=v^{q},\;\;\;\;\mbox{in}\;\mathbb{R}^{n-1}.
$$
 in \cite{GNW}, \cite{L} to the square root of Laplacian operator case.

In the general case $a\in (-1,0)$, Fang, Gui and the author in \cite{Du1} obtain existence of positive  axial symmetric  solutions  in
space dimensions $n\geq3$  for
$
 q\geq\frac{n-a}{n+a-2}
$.
Further the  author   establish a unique singular positive axial symmetric  solution $\psi_\infty(x)$,
for   $q\geq\frac{n-1}{n+a-2}$,
in $n\geq3$ for all $a\in (-1,0)\cup(0,1)$  (see \cite{Du2}).

Naturally we hope to establish  similar qualitative properties
of the positive solution of (\ref{e2}) for the general case $a\in (-1,0)$
as that for the case $a=0$ in \cite{H}.
Our main results in this paper are stated as follows.

\begin{thm}\label{thm1}
Let $q\geq\frac{n-a}{n+a-2}$ and denote $m_q=\frac{1-a}{q-1}$. Then there exists a family of positive axial symmetric solutions $u_\beta(x)$ ($u_\beta(0)=\beta>0$) of (\ref{e2}) satisfying the following properties

(i)$
u_\beta(x)=\beta u_1(\beta^{\frac{q-1}{1-a}}x);
$

(ii)$
u_\beta(x)\leq C(1+|x|)^{-m_q};
$

(iii)$
\lim_{\beta\rightarrow\infty}u_\beta(x)=\psi_\infty(x)
$ for $x\in \mathbb{R}_{+}^{n}\backslash \{0\}$, $q\neq\frac{n-a}{n+a-2}$.
\end{thm}

The following Theorem is the asymptotic expansion for the JL-supercritical,   JL-critical
and JL-subcritical case (see Definition \ref{defi1}).

\begin{thm}\label{thm2}
(I) Let $q$ be JL-supercritical or JL-critical. Then there exists $\varepsilon>0$ such that for any
$\beta>0$ there exist $C_1(\beta)<0, C_2(\beta)\in\mathbb{R}$, such that the solutions $u_\beta$ of (\ref{e2}) hold the asymptotic expansion for large $r>0$
$$
u_\beta(x)=\psi_\infty(x)+[C_1+\mbox{O}(r^{-\varepsilon})]e_1(\theta)r^{-(m_q+\rho_1)}, ~~~~\mbox{for JL-supercritical case},
$$
and
$$
u_\beta(x)=\psi_\infty(x)+[C_1 \ln r+C_2+\mbox{O}(r^{-\varepsilon})]e_1(\theta)r^{-\frac{n+a-2}{2}}, ~~~~\mbox{for JL-critical case},
$$
where we have used the polar coordinate $|x'|=r\sin \theta, ~x_n=r\cos \theta$, and
$$
\rho_1=\frac{(n+a-2-2m_q)-\sqrt{(n+a-2-2m_q)^2+4[m_q(n+a-2-m_q)+\lambda_1]}}{2},
$$
where $\lambda_i,~e_i(\theta)$ are the $i$-th eigenvalue, the $i$-th eigenfunction  of (\ref{e23})
with
$$\int_0^{\frac{\pi}{2}}e_i(\theta)e_j(\theta)\sin^{n-2}\cos^a\theta d\theta=\delta_{ij},~~~~e_i(\frac{\pi}{2})>0.$$
Moreover the asymptotic expansion holds uniformly for $\theta\in(0,\frac{\pi}{2}).$

(II) Let $q$ be JL-subcritical. Then the solutions $u_\beta$ of (\ref{e2}) satisfies one of the following two asymptotic expansions for large $r>0$

1) there exists $(c_1,c_2)\neq0$ and $\varepsilon>0$ such that
$$
u_\beta(x)=\psi_\infty(x)+[c_1 \sin(K\ln r)+c_2\cos(K\ln r)+\mbox{O}(r^{-\varepsilon})]e_1(\theta)r^{-\frac{n+a-2}{2}};
$$

2) there exists $c>0$ and $\varepsilon>0$ such that
$$
u_\beta(x)=\psi_\infty(x)+[c+\mbox{O}(r^{-\varepsilon})]e_2(\theta)r^{-(m_q+\rho_2)},
$$
where $K>0, ~\rho_2>0$ are given by
$$
K=\frac{\sqrt{-(n+a-2-2m_q)^2-4[m_q(n+a-2-m_q)+\lambda_1]}}{2},
$$
$$
\rho_2=\frac{(n+a-2-2m_q)+\sqrt{(n+a-2-2m_q)^2+4[m_q(n+a-2-m_q)+\lambda_2]}}{2}.
$$
Moreover the asymptotic expansion holds uniformly for $\theta\in(0,\frac{\pi}{2}).$
\end{thm}

Therefore,  by using the Caffarelli-Silvestre extension, we know that  our results generalize the results from the square root Laplacian operator to general fractional Laplacian operator
$$
(-\Delta)^{s}v=v^{q},\;\;\;\;\mbox{in}\;\mathbb{R}^{n-1},
$$
where $s=\frac{1-a}{2}\in (\frac{1}{2},1)$, since $a\in (-1,0)$.

The paper is organized as follows. Section 2 contains some notations and an estimate of the infimum $C_a$ in generalized trace
Hardy inequality, which will be used in the definition of
  JL-critical exponent in Section 3.
In section 4 estimates of eigenvalues of an eigenvalue problem are made, which will be used in  asymptotic expansions for positive axial symmetric  solutions.
In section 5 we derive decay estimates and limit behavior of them.
Section 6 contains the proof of  asymptotic expansions for them.

\section{Preliminaries}

We say that a function $u$ is axially symmetric with respect to the $x_n$-axis, if it can be expressed by $u(x)=\tilde{u}(|x'|,x_n)$ for
some function $\tilde{u}$. For axially symmetric function, one need to introduce the polar coordinates $(r, \theta)$
 \begin{equation}\label{e4}
 |x'|=r\sin \theta, ~~~~~x_n=r\cos \theta,~~~\theta\in [0,\frac{\pi}{2}].
\end{equation}

We first recall the main results in \cite{Du1} and \cite{Du2}.

\begin{pro} \label{prop1}(\cite{Du1})
Let $q\geq\frac{n-a}{n+a-2}$ and $a\in(-1,0)$. The  problem (\ref{e2}) admits a positive axially symmetric solution  $u(x)$ satisfying
 $u_\theta>0$ and $u(0)=1$.
\end{pro}

In \cite{Du2}, the author consider the existence of positive singular solutions of the form
$$
\psi_\infty(x)=V(\theta)r^{\frac{a-1}{q-1}}.
$$
To obtain $V(\theta)$, one need to solve
\begin{eqnarray}\label{e5}
\left\{\begin{array}{l}(\sin^{n-2} \theta\cos^a \theta V_\theta)_\theta=\gamma V\sin^{n-2} \theta\cos^a \theta\;\;\;\;\;\theta\in (0,\frac{\pi}{2}),\\
\lim_{\theta\rightarrow \frac{\pi}{2}}V_\theta(\theta)\cos^a \theta =V^{q}(\frac{\pi}{2}),\\
\end{array}
\right.
\end{eqnarray}
where
$$
\gamma:=m_q(n+a-2-m_q)>0.
$$

\begin{pro} \label{prop2}(\cite{Du2})
Let $q\geq\frac{n-1}{n+a-2}$ and $a\in(-1,0)\cup (0,1)$. The  problem (\ref{e5}) admits a unique solution  $V(\theta)$ satisfying
 $V_\theta(\theta)>0$.
\end{pro}

Next we introduce some notations that will be used.
Let $S^{n-1}_+=\{\omega=(\omega',\cos\theta): \omega'\in\mathbb{R}^{n-1}, \theta\in[0,\frac{\pi}{2}), |\omega'|^2+\cos^2\theta=1\}$
be a half unit sphere. For $p>1$, we define the space
$$L^p_{\mbox{sym}}(S^{n-1}_+)=\left\{e: e(\omega)~ \mbox{depends only on}~ \theta~ \mbox{and}
\int_0^{\frac{\pi}{2}}|e(\theta)|^p\sin^{n-2} \theta\cos^a \theta d\theta<+\infty \right\}.$$
We denote
$$
\|e(\theta)\|_{p,a}:=\left(\int_0^{\frac{\pi}{2}}|e(\theta)|^p\sin^{n-2} \theta\cos^a \theta d\theta\right)^{\frac{1}{p}}.
$$
Further we define the space
$$\mathcal{H}^k_a(S^{n-1}_+)=\left\{e\in L^p_{\mbox{sym}}(S^{n-1}_+):
\|e\|_{\mathcal{H}^k_a(S^{n-1}_+)}:=\sum_{i=0}^k\|\partial^i_\theta e(\theta)\|_{p,a}<+\infty \right\}.$$

We also define the following operator on $S^{n-1}_+$ according to the norm $\|\cdot\|_{p,a}$
$$
\Delta_{S,a}:=\partial_{\theta\theta}+[(n-2)\cot\theta-a\tan\theta]\partial_\theta.
$$
Note that
$$
\sin^{n-2} \theta\cos^a \theta\Delta_{S,a}e(\theta)=(\sin^{n-2} \theta\cos^a \theta e_\theta)_\theta.
$$

We will establish the following elementary lemma.

\begin{lem}\label{lemm1}
For $a\in(-1,0)$, we have
\begin{equation}\label{e6}
\int_0^{\frac{\pi}{2}}\sin^{n+a-2} \theta d\theta<\int_0^{\frac{\pi}{2}}\sin^{n-2} \theta\cos^a \theta d\theta.
\end{equation}
\end{lem}
\emph{Proof.} Since
$$
\int_0^{\frac{\pi}{2}}\sin^{n+a-2} \theta d\theta=\int_0^{\frac{\pi}{4}}\sin^{n+a-2} \theta d\theta
+\int_0^{\frac{\pi}{4}}\cos^{n+a-2} \theta d\theta
$$
and
$$
\int_0^{\frac{\pi}{2}}\sin^{n-2} \theta \cos^a \theta d\theta=\int_0^{\frac{\pi}{4}}\sin^{n-2} \theta \cos^a \theta d\theta
+\int_0^{\frac{\pi}{4}}\cos^{n-2} \theta\sin^a \theta d\theta,
$$
we have
\begin{eqnarray}
\nonumber
&&\int_0^{\frac{\pi}{2}}\sin^{n+a-2} \theta d\theta-\int_0^{\frac{\pi}{2}}\sin^{n-2} \theta\cos^a \theta d\theta\\
\nonumber
&=& \int_0^{\frac{\pi}{4}}(\sin^{n-2} \theta-\cos^{n-2} \theta) (\sin^a \theta-\cos^a \theta)d\theta<0.
\end{eqnarray}

The lemma is proved.

We introduce the space
$
H^1_a(\mathbb{R}_{+}^{n}):=\left\{u: \int_{\mathbb{R}_{+}^{n}}(u^2+|\nabla u|^2)x_{n}^{a}dx<+\infty\right\},
$
and consider the following minimizing problem
$$
C_a:=\inf_{u\in H^1_a(\mathbb{R}_{+}^{n})}\frac{\int_{\mathbb{R}_{+}^{n}}|\nabla u|^2x_{n}^{a}dx}{\int_{\partial\mathbb{R}_{+}^{n}}| u|^2|x'|^{a-1}dx'},
$$
which plays an essential role to define a so-called JL-critical exponent in Section 3.

\begin{lem}\label{lemm2}
Let $a\in(-1,0)$ and $n\geq4$. We have
$$
C_a>\frac{n+a-3}{2}.
$$
\end{lem}
\emph{Proof.} One may assume that $u=u(\bar{r},t)$ where $\bar{r}=|x'|=|(x_1,\cdots, x_{n-1})|$ and $t=x_n$.
Then
$$
\int_{\partial\mathbb{R}_{+}^{n}}| u|^2|x'|^{a-1}dx'=(n-1)\omega_{n-1}\int_0^\infty u^2(\bar{r},0)\bar{r}^{n+a-3}d\bar{r},
$$
where $\omega_{n-1}$ denotes the measure of the unit sphere in $\mathbb{R}^{n-1}$.
Since
$$
u^2(\bar{r},0)=-2\int_0^\infty u(\bar{r},t)\partial_t u(\bar{r},t)dt,
$$
Hence
\begin{eqnarray}
\nonumber
&&\int_{\partial\mathbb{R}_{+}^{n}}| u|^2|x'|^{a-1}dx'=-2(n-1)\omega_{n-1}\int_0^\infty\int_0^\infty u(\bar{r},t)\partial_t u(\bar{r},t)
\bar{r}^{n+a-3}d\bar{r}dt\\
\nonumber
&\leq& 2(n-1)\omega_{n-1}\int_0^\infty\left(\int_0^\infty u^2(\bar{r},t)\bar{r}^{n+a-4}d\bar{r}\right)^{1/2}
\left(\int_0^\infty (\partial_t u(\bar{r},t))^2\bar{r}^{n+a-2}d\bar{r}\right)^{1/2}dt.
\end{eqnarray}
We use now the inequality
$$
\int_0^\infty u^2(\bar{r},t)\bar{r}^{n+a-4}d\bar{r}\leq\frac{4}{(n+a-3)^2}\int_0^\infty (\partial_{\bar{r}} u(\bar{r},t))^2\bar{r}^{n+a-2}d\bar{r},
$$
which follows from the generalized Hardy's inequality
$$
\int_{\mathbb{R}^{d}}\frac{ u^2}{|x|^{2(b+1)}}dx\leq\frac{4}{(d-2-2b)^2}\int_{\mathbb{R}^{d}}\frac{|\nabla u|^2}{|x|^{2b}}dx.
$$
We obtain
\begin{eqnarray}
\nonumber
&&\int_{\partial\mathbb{R}_{+}^{n}}| u|^2|x'|^{a-1}dx'\\
\nonumber
&\leq& \frac{4(n-1)\omega_{n-1}}{n+a-3}\int_0^\infty\left(\int_0^\infty (\partial_{\bar{r}} u(\bar{r},t))^2\bar{r}^{n+a-2}d\bar{r}\right)^{1/2}
\left(\int_0^\infty (\partial_t u(\bar{r},t))^2\bar{r}^{n+a-2}d\bar{r}\right)^{1/2}dt\\
\label{e7}
&\leq&\frac{2(n-1)\omega_{n-1}}{n+a-3}\int_0^\infty\int_0^\infty [(\partial_{\bar{r}} u(\bar{r},t))^2+(\partial_t u(\bar{r},t))^2]\bar{r}^{n+a-2}d\bar{r}dt.
\end{eqnarray}
We introduce the coordinates $(r,\theta)$
$$
\bar{r}=r\sin\theta,~~~~~~t=r\cos\theta.
$$
Then one has
\begin{eqnarray}
\nonumber
&&
\int_0^\infty\int_0^\infty [(\partial_{\bar{r}} u(\bar{r},t))^2+(\partial_t u(\bar{r},t))^2]\bar{r}^{n+a-2}d\bar{r}dt
\\
\label{e8}
&=&
\int_0^\infty r^{n+a-1}\left(\int_0^{\frac{\pi}{2}} |\nabla u(r\sin\theta,r\cos\theta)|^2\sin^{n+a-2}d\theta\right) dr.
\end{eqnarray}
We may further assume that $|\nabla u(r\sin\theta,r\cos\theta)|^2=|\nabla u(r\cos\theta,r\sin\theta)|^2$. From this, using
the similar argument of Lemma \ref{lemm1}, we obtain
\begin{eqnarray}
\nonumber
&&
\int_0^\infty r^{n+a-1}\left(\int_0^{\frac{\pi}{2}} |\nabla u(r\sin\theta,r\cos\theta)|^2\sin^{n+a-2}d\theta\right) dr
\\
\label{e9}
&<&
\int_0^\infty r^{n+a-1}\left(\int_0^{\frac{\pi}{2}} |\nabla u(r\sin\theta,r\cos\theta)|^2\sin^{n-2}\cos^a\theta d\theta\right) dr.
\end{eqnarray}
Note that
\begin{eqnarray}
\nonumber
&&
\int_0^\infty r^{n+a-1}\left(\int_0^{\frac{\pi}{2}} |\nabla u(r\sin\theta,r\cos\theta)|^2\sin^{n-2}\cos^a\theta d\theta\right) dr
\\ \label{e10}
&&=
\frac{1}{(n-1)\omega_{n-1}}\int_{\mathbb{R}_{+}^{n}}|\nabla u|^2x_{n}^{a}dx.
\end{eqnarray}
Hence, from (\ref{e7})-(\ref{e10}), we have
$$
\int_{\partial\mathbb{R}_{+}^{n}}| u|^2|x'|^{a-1}dx'\leq \frac{2}{n+a-3}\int_{\mathbb{R}_{+}^{n}}|\nabla u|^2x_{n}^{a}dx.
$$

The proof of this lemma is complete.

\begin{mrem}\label{mrem0}
This lemma   generalize the conclusion that $C_0>\frac{n-3}{2}$ for the case $a=0$, which is obtained by D\'{a}vila-Dupaigne-Montenegro in \cite{Da}.
\end{mrem}

We define
$$h_{n,a}:=-\frac{(n+a-2)^2}{4}.$$
Next we show that $h_{n,a}$ is the first eigenvalue of the following eigenvalue problem
\begin{equation}\label{e11}
-\Delta_{S,a}e=\lambda e, ~~~\theta\in(0,\frac{\pi}{2}),~~~\lim_{\theta\rightarrow \frac{\pi}{2}}e_\theta(\theta)\cos^a \theta =C_ae(\frac{\pi}{2}).
\end{equation}
To this end, we introduce the space
$$
H^1_a(S^{n-1}_+):=\left\{e: \int_{ S^{n-1}_+}(e^2+|\nabla_{S} e|^2)d\omega=\int_{ S^{n-1}_+}(e^2+|\nabla_{S} e|^2)\sin^{n-2}\cos^a\theta d\theta d\omega'<+\infty\right\},
$$
where $\nabla_{S}e$ is a gradient of $e$ on $S^{n-1}_+$. We consider a minimizing problem
\begin{equation}\label{e12}
C_S:=\inf_{e\in H^1_a(S^{n-1}_+)}\frac{\|\nabla_{S} e\|^2_{L^2_a(S^{n-1}_+)}-h_{n,a}\| e\|^2_{L^2_a(S^{n-1}_+)}}{\| e\|^2_{L^2(\partial S^{n-1}_+)}},
\end{equation}
where $\| e\|^2_{L^2_a(S^{n-1}_+)}:=\int_{ S^{n-1}_+}e^2(\theta)\sin^{n-2}\cos^a\theta d\theta d\omega'$.
For simplicity of notations, hereafter we set
$$
d\mu:=\sin^{n-2}\cos^a\theta d\theta,~~~e_B:=e(\frac{\pi}{2}).
$$
By a trace inequality, it is verified that there exists a minimizer $e\in\mathcal{H}^1_a(S^{n-1}_+)$
of (\ref{e12}) satisfying $e(\theta)>0$ and
\begin{equation}\label{e13}
-\Delta_{S,a}e=h_{n,a} e, ~~~\theta\in(0,\frac{\pi}{2}),~~~\lim_{\theta\rightarrow \frac{\pi}{2}}e_\theta(\theta)\cos^a \theta =C_Se_B.
\end{equation}
Since $e(\theta)>0$, then $h_{n,a}$ is the first eigenvalue of this problem. Therefore, for problem  (\ref{e11}), to prove $\lambda_1=h_{n,a}$  we only need to show $C_S=C_a$.

\begin{lem}\label{lemm3}
 $C_S=C_a$.
\end{lem}
\emph{Proof.}
Let $u\in C_c^1(\overline{\mathbb{R}_{+}^{n}})$. We have
$$
\int_{\partial\mathbb{R}_{+}^{n}}| u(x',0)|^2|x'|^{a-1}dx'=\int_0^\infty \|u(r,\cdot)\|^2_{L^2(\partial S^{n-1}_+)}r^{n+a-3}dr.
$$
By the definition of $C_S$ and the generalized Hardy inequality, one has
\begin{eqnarray}
\nonumber
C_S\int_0^\infty \|u(r,\cdot)\|^2_{L^2(\partial S^{n-1}_+)}r^{n+a-3}dr&\leq & \int_0^\infty\int_{ S^{n-1}_+} (|\nabla_{S}u|^2-
h_{n,a}u^2)r^{n+a-3}d\omega dr
\\
\nonumber
&\leq&
\int_0^\infty\int_{ S^{n-1}_+} (r^{-2}|\nabla_{S}u|^2+
u_r^2)r^{n+a-1}d\omega dr.
\end{eqnarray}
Note that $|\nabla u|^2=u_r^2+r^{-2}|\nabla_{S}u|^2$, we obtain
$$
C_S\int_0^\infty \|u(r,\cdot)\|^2_{L^2(\partial S^{n-1}_+)}r^{n+a-3}dr\leq
\int_{\mathbb{R}_{+}^{n}}|\nabla u|^2x_n^adx,
$$
which gives that $C_a\geq C_S$.

Next we show $C_a\leq C_S$. We need to construct a sequence $\{u_i\}_{i\in N}\subset H^1_a(\mathbb{R}_{+}^{n})$ such that
\begin{equation}\label{e14}
\lim_{i\rightarrow\infty}\frac{\int_{\mathbb{R}_{+}^{n}}|\nabla u_i|^2x_{n}^{a}dx}{\int_{\partial\mathbb{R}_{+}^{n}}| u_i|^2|x'|^{a-1}dx'}=C_S.
\end{equation}
To this end, we set
\begin{equation}\label{e15}
u_i(x)=\left\{\begin{array}{l}e(\theta)i^{(n+a-2)/2}  ~~~~~~~~~~~~ \mbox{if}~~
r\in [0,1/i),\\
e(\theta)r^{(n+a-2)/2}  ~~~~~~~~~~~ \mbox{if}~~
r\in [1/i,1),\\
\chi(r)e(\theta)r^{(n+a-2)/2}  ~~~~~~ \mbox{if}~~
r\in [1,\infty),\\
\end{array}
\right.
\end{equation}
where $e(\theta)$ is a minimizer of (\ref{e12}) and $\chi(r)$ is a cut-off function such that $\chi(r)=1$
for $r\in[0,1]$ and $\chi(r)=0$ for $r\in[2,\infty)$.
Calculation shows that
$$
\int_{\mathbb{R}_{+}^{n}}|\nabla u_i|^2x_{n}^{a}dx=C_S\|e\|^2_{L^2(\partial S^{n-1}_+)}(\ln i)+\mbox{O}(1)
$$
and
$$
\int_{\partial\mathbb{R}_{+}^{n}}| u(x',0)|^2|x'|^{a-1}dx'=\|e\|^2_{L^2(\partial S^{n-1}_+)}(\ln i)+\mbox{O}(1).
$$
Hence we obtain (\ref{e14}). We complete the proof of this lemma.

\section{JL-critical exponent}

We introduce
\begin{equation}\label{e16}
J(q):=\inf_{u\in H^1_a(\mathbb{R}_{+}^{n})}\frac{\int_{\mathbb{R}_{+}^{n}}|\nabla u|^2x_{n}^{a}dx-q\int_{\partial\mathbb{R}_{+}^{n}}\psi_\infty^{q-1}u^2dx'}{\int_{\partial\mathbb{R}_{+}^{n}}| u|^2|x'|^{a-1}dx'}.
\end{equation}
Substituting the explicit expression of $\psi_\infty(x)=V(\theta)r^{\frac{a-1}{q-1}}$ into (\ref{e16}), we have
$$
J(q)=C_a-qV_B^{q-1}.
$$

\begin{defi}\label{defi1}
We say that an exponent $q$ is JL-supercritical if $J(q)>0$,   JL-critical if $J(q)=0$ and JL-subcritical if $J(q)<0$.
\end{defi}

\begin{mrem}\label{mrem1}
Hence, an exponent $q$ is JL-supercritical if $C_a>qV_B^{q-1}$,   JL-critical if $C_a=qV_B^{q-1}$ and JL-subcritical if $C_a<qV_B^{q-1}$.
\end{mrem}

Next we will show that $q$ is JL-supercritical if $q$ and $n$ are large enough, and $q$ is JL-subcritical if $q$ is close to $\frac{n-a}{n+a-2}$.

\begin{lem}\label{lemm4}
There exists  $n_0\in N$ and $q_1>\frac{n_0-a}{n_0+a-2}$ such that for $n\geq n_0$ and $q>q_1$ we have $C_a>qV_B^{q-1}$.
\end{lem}
\emph{Proof.} Integrating (\ref{e5}) over $(0,\frac{\pi}{2})$, we have
\begin{equation}\label{e17}
V_B^q=\lim_{\theta\rightarrow \frac{\pi}{2}}V_\theta(\theta)\cos^a \theta =\gamma\int_0^{\frac{\pi}{2}}V(\theta)d\mu\leq
\gamma\int_0^{\frac{\pi}{2}}V_Bd\mu,
\end{equation}
where in the last inequality we used the fact that $V_\theta>0$.
Hence we have
$$
qV_B^{q-1}\leq q\gamma\int_0^{\frac{\pi}{2}}d\mu=qm_q(n+a-2-m_q)I_{n,a},
$$
where $I_{n,a}:=\int_0^{\frac{\pi}{2}}d\mu=\int_0^{\frac{\pi}{2}}\sin^{n-2}\cos^a\theta d\theta\rightarrow0$ as $n\rightarrow+\infty$.
Note that for $q>q_1$ we have $qm_q=q\frac{1-a}{q-1}\leq C$, where $C$ is a positive constant independence of $n$.
On the other hand, by Lemma \ref{lemm2}, one has $C_a>\frac{n+a-3}{2}$.
So $qm_q(n+a-2-m_q)I_{n,a}< C_a$ for large $n$.
 Hence the result of this lemma is true.

\begin{lem}\label{lemm5}
For $n\geq3$, there exists   $q_0\geq\frac{n-a}{n+a-2}$ such that  $C_a<qV_B^{q-1}$ for $q\in[\frac{n-a}{n+a-2},q_0)$.
Furthermore, in these low dimensions $n=3,4,5,6$, we have $C_a<qV_B^{q-1}$ for any $q\geq\frac{n-a}{n+a-2}$.
\end{lem}
\emph{Proof.} By the explicit expression of $\gamma$ and $h_{n,a}$,
it is easily seen that $h_{n,a}\geq\gamma$ for $q\in[\frac{n-a}{n+a-2},+\infty)$, moreover, the strict inequality holds unless $q=\frac{n-a}{n+a-2}$.
Let $e(\theta)$ be a positive solution of (\ref{e13}) with $e_B=V_B$. We set $W(\theta):=e(\theta)-V(\theta)$, then
we have
\begin{equation}\label{e18}
\left\{\begin{array}{l}W_{\theta\theta}+(n-2)\cot\theta W_\theta-a\tan\theta W_\theta\geq\gamma W\;\;\;\;\;\mbox{in}\;\;(0,\frac{\pi}{2}),\\
W=0 \;\;\;\;\;\;\;\;
\;\;\;\;\;\;\;\;\;~~~~~~~~~~~~~~~~~~~~~~~~~~~~~~~~~~~\mbox{on}\;\;\{\frac{\pi}{2}\}.\\
\end{array}
\right.
\end{equation}

We claim that
\begin{equation}\label{e19}
W(\theta)\leq0~~~~~~\mbox{in}~~(0,\frac{\pi}{2}).
\end{equation}
First we suppose that $W(0)>0$.
Note that $\lim_{\theta\rightarrow0}\frac{W_\theta(\theta)}{\tan\theta}=\lim_{\theta\rightarrow0}\frac{W_\theta(\theta)-W_\theta(0)}{\theta}=W_{\theta\theta}(0)$.
By this and (\ref{e18}), we have $(n-1)W_{\theta\theta}(0)\geq\gamma W(0)>0$.
Since $W_\theta(0)=0$, there exists $\delta>0$ such that $W_\theta>0$ in $(0,\delta)$. Due to $W(\frac{\pi}{2})=0$,
there exists $\theta_0\in(0,\frac{\pi}{2})$ such that $W(\theta_0)>0,~W_\theta(\theta_0)=0$ and $W_{\theta\theta}(\theta_0)\leq0$.
However, these contradict with (\ref{e18}). Hence we obtain $W(0)\leq0$. Now we suppose (\ref{e19}) is false. Then there exists $\theta_1\in(0,\frac{\pi}{2})$
such that $W(\theta_1)>0,~W_\theta(\theta_1)=0$ and $W_{\theta\theta}(\theta_1)\leq0$,
which contradict with (\ref{e18}) again. Hence the claim (\ref{e19}) is proved.

Multiplying the equation in the following problem
$$
-\Delta_{S,a}e=h_{n,a} e, ~~~\theta\in(0,\frac{\pi}{2}),~~~\lim_{\theta\rightarrow \frac{\pi}{2}}e_\theta(\theta)\cos^a \theta =C_ae_B,
$$
by $\sin^{n-2}\cos^a\theta$ and integrating  over $(0,\frac{\pi}{2})$,
we obtain
\begin{equation}\label{e20}
C_ae_B=h_{n,a}\int_0^{\frac{\pi}{2}}e(\theta)d\mu.
\end{equation}
From  (\ref{e17}), (\ref{e20}) and the assumption $e_B=V_B$, we have
$$
qV_B^{q-1}-C_a=V_B^{-1}(qV_B^{q}-C_ae_B)=V_B^{-1}\left(q\gamma\int_0^{\frac{\pi}{2}}V(\theta)d\mu-h_{n,a}\int_0^{\frac{\pi}{2}}e(\theta)d\mu\right).
$$
From (\ref{e19}) we have
$$
qV_B^{q-1}-C_a\geq V_B^{-1}(q\gamma-h_{n,a})\int_0^{\frac{\pi}{2}}e(\theta)d\mu.
$$
Then to prove $qV_B^{q-1}>C_a$, we need to show $q\gamma> h_{n,a}$. From the explicit expression of $\gamma$, we have
\begin{equation}\label{e21}
q\gamma-h_{n,a}=\left(1-a+\frac{1-a}{q-1}\right)\left((n+a-2)-\frac{1-a}{q-1}\right)-\left(\frac{n+a-2}{2}\right)^2.
\end{equation}
Now we choose $q=\frac{n-a}{n+a-2}$ in (\ref{e21}), then we see that $q\gamma-h_{n,a}>0$.
Hence by a continuity principle, there exists $q_0>\frac{n-a}{n+a-2}$ such that $q\gamma-h_{n,a}>0$ for $q\in[\frac{n-a}{n+a-2},q_0)$,
which shows the first statement. 

Next we prove the second statement. We regard the right-hand side of (\ref{e21})
as a function $G$ of $\frac{1-a}{q-1}$. Then $q\geq\frac{n-a}{n+a-2}$ is equivalent to $\frac{1-a}{q-1}\in(0,\frac{n+a-2}{2}]$.
Hence
\begin{equation}\label{e22}
\inf_{q\geq\frac{n-a}{n+a-2}}(q\gamma-h_{n,a})=\inf_{\tau\in(0,\frac{n+a-2}{2}]}G(\tau),
\end{equation}
where $G(\tau):=(1-a+\tau)((n+a-2)-\tau)-\frac{(n+a-2)^2}{4}$.
Elementary computation shows that, for $n=3$
$$
\inf_{\tau\in(0,\frac{n+a-2}{2}]}G(\tau)=G(\frac{n+a-2}{2})=\frac{1-a^2}{2}>0,
$$
and for $n=4$
$$
\inf_{\tau\in(0,\frac{n+a-2}{2}]}G(\tau)=G(\frac{n+a-2}{2})=\frac{(2+a)(1-a)}{2}>0.
$$
For the case $n=5$, if $a\in(-1,-\frac{1}{3}]$, we have
$$
\inf_{\tau\in(0,\frac{n+a-2}{2}]}G(\tau)=G(\frac{n+a-2}{2})=\frac{(3+a)(1-a)}{2}>0,
$$
and if $a\in (-\frac{1}{3},0)$, we have
$$
\inf_{\tau\in(0,\frac{n+a-2}{2}]}G(\tau)=G(0)=\frac{(3+a)(1-a)}{4}>0.
$$
For the case $n=6$, if $a\in(-1,-\frac{2}{3}]$, we have
$$
\inf_{\tau\in(0,\frac{n+a-2}{2}]}G(\tau)=G(\frac{n+a-2}{2})=\frac{(4+a)(1-a)}{2}>0,
$$
and if $a\in (-\frac{2}{3},0)$, we have
$$
\inf_{\tau\in(0,\frac{n+a-2}{2}]}G(\tau)=G(0)=\frac{(4+a)(-a)}{4}>0.
$$
As for the case $n\geq7$, we have $\inf_{\tau\in(0,\frac{n+a-2}{2}]}G(\tau)=G(0)<0.$

The proof of this lemma is complete.

\section{Estimate of eigenvalues}

In  this section we will make some estimates for the eigenvalues of the following problem.
\begin{equation}\label{e23}
-\Delta_{S,a}e=\lambda e, ~~~\theta\in(0,\frac{\pi}{2}),~~~\lim_{\theta\rightarrow \frac{\pi}{2}}e_\theta(\theta)\cos^a \theta =qV_B^{q-1}e_B.
\end{equation}
We denote by $\lambda_i,~e_i(\theta)$ the $i$-th normalized eigenfunction such that
$$
\int_0^{\frac{\pi}{2}}e_i(\theta)e_j(\theta)d\mu=\delta_{ij},~~~e_{iB}>0.
$$

First we give estimates for the first eigenvalue for this problem.
\begin{lem}\label{lemm6}
 $\lambda_1<-\gamma$ and

 $\lambda_1>-\frac{(n+a-2)^2}{4}$~~~~if $q$ is JL-supercritical;

  $\lambda_1=-\frac{(n+a-2)^2}{4}$~~~~if $q$ is JL-critical;

   $\lambda_1<-\frac{(n+a-2)^2}{4}$~~~~if $q$ is JL-subcritical.
\end{lem}
\emph{Proof.} For the first eigenvalue $\lambda_1$   of (\ref{e23}), it is characterized by
\begin{equation}\label{e24}
\lambda_1=\inf_{e\in \mathcal{H}^1_a(S^{n-1}_+)}\frac{\|\partial_{\theta} e\|^2_{2,a}-qV_B^{q-1} e^2_B}{\| e\|^2_{2,a}}.
\end{equation}
Recall that $V(\theta)$ is a positive solution of (\ref{e5}), then $-\gamma$ is the first eigenvalue of the problem
$$
-\Delta_{S,a}e=\lambda e, ~~~\theta\in(0,\frac{\pi}{2}),~~~\lim_{\theta\rightarrow \frac{\pi}{2}}e_\theta(\theta)\cos^a \theta =V_B^{q-1}e_B.
$$
Hence
\begin{equation}\label{e25}
-\gamma=\inf_{e\in \mathcal{H}^1_a(S^{n-1}_+)}\frac{\|\partial_{\theta} e\|^2_{2,a}-V_B^{q-1} e^2_B}{\| e\|^2_{2,a}}.
\end{equation}
By   (\ref{e24}), (\ref{e25}) and the fact $q>1$, we obtain that $\lambda_1<-\gamma$.

In  Section 2 we have proved that $h_{n,a}=-\frac{(n+a-2)^2}{4}$ is the first eigenvalue of (\ref{e11}),
hence  we have
\begin{equation}\label{e26}
-\frac{(n+a-2)^2}{4}=\inf_{e\in \mathcal{H}^1_a(S^{n-1}_+)}\frac{\|\partial_{\theta} e\|^2_{2,a}-C_a e^2_B}{\| e\|^2_{2,a}}.
\end{equation}
By Remark \ref{mrem1},  (\ref{e24}) and (\ref{e26}), we obtain the rest results of this lemma.

We set
$$\sigma:=(n+a-2)-2m_q.$$
 Elementary computation shows that $\sigma^2+4(\gamma+\lambda_1)=(n+a-2)^2+4\lambda_1$.
From this and Lemma \ref{lemm6}, we obtain the following lemma.

\begin{lem}\label{lemm7}
 One has

 $\sigma^2+4(\gamma+\lambda_1)>0$~~~~if $q$ is JL-supercritical;

  $\sigma^2+4(\gamma+\lambda_1)=0$~~~~if $q$ is JL-critical;

   $\sigma^2+4(\gamma+\lambda_1)<0$~~~~if $q$ is JL-subcritical.
\end{lem}

Next we show that the second eigenvalue of (\ref{e23}) is positive.

\begin{lem}\label{lemm8}
It holds that $\lambda_2>0.$
\end{lem}
\emph{Proof.} Let $e_2(\theta)$ be a a corresponding eigenfunction  with $e_2(0)<0$.
By Strum's comparison theorem, $e_2(\theta)$ has just one zero in $(0,\frac{\pi}{2})$, which is denoted by $\theta_0$.
Then $\partial_\theta e_2(\theta_0)\geq0$. Clearly $\partial_\theta e_2(\theta_0)\neq0$. Hence
$\partial_\theta e_2(\theta_0)>0$.
Multiplying (\ref{e23}) by $\sin^{n-2}\theta_0\cos^a\theta_0$ and integrating over $(0,\theta_0)$, we obtain
$$
\partial_\theta e_2(\theta_0)=\frac{-\lambda_2}{\sin^{n-2}\theta_0\cos^a\theta_0}\int_0^{\theta_0}e_2(\theta)d\mu,
$$
which gives that $\lambda_2>0.$

\section{Decay estimates and Limit behavior}

In this section we will prove Theorem \ref{thm1}.

Recall that from Proposition \ref{prop1}, (\ref{e2}) admits a positive axial symmetric solution $u(x)$
satisfying $u(0)=1$. We set $u_1(x):=u(x)$ and
$$
u_\beta(x):=\beta u_1(\beta^{\frac{q-1}{1-a}}x).
$$
Then $u_\beta(x)$ is a positive axial symmetric solution of (\ref{e2})
satisfying $u_\beta(0)=\beta$. Hence Theorem \ref{thm1} (i) is proved.

Now we will  follow the idea in \cite{H} to prove Theorem \ref{thm1} (ii).
We denote a positive axial symmetric solution of (\ref{e2}) as $u(x)=u(r,\theta)$, which satisfies $u_\theta>0$. Then we need to prove the decay estimate
\begin{equation}\label{e27}
u(r,\theta)=u(x)\leq C(1+|x|)^{-m_q}.
\end{equation}

We define
$$
v(t,\theta):=r^{\frac{1-a}{q-1}}u(r,\theta),~~~~r:=e^t.
$$
Then $v(t,\theta)$ satisfies
\begin{equation}\label{e28}
\left\{\begin{array}{l}v_{tt}+\sigma v_t-\gamma v+\Delta_{S,a}v=0\;\;\;\;\;\mbox{in}\;\;\mathbb{R}\times (0,\frac{\pi}{2}),\\
\lim_{\theta\rightarrow \frac{\pi}{2}}v_\theta\cos^a \theta =v_B^{q}(t) \;\;\;\;\;\;\;\;
\;\mbox{on}\;\; \mathbb{R},\\
\end{array}
\right.
\end{equation}
where $\sigma$ is the positive constant defined in the previous section.
To prove the decay estimate (\ref{e27}), we need to show that $v(t,\theta)$ is bounded on $\mathbb{R}_+\times (0,\frac{\pi}{2})$.
We set
$$
\tilde{v}(t):=\int_0^{\frac{\pi}{2}}v(t,\theta)d\mu.
$$
We will first show that $\tilde{v}(t)$ is bounded on $\mathbb{R}_+$.

Multiplying (\ref{e28}) by $\sin^{n-2}\theta\cos^a\theta$ and integrating with respect to $\theta$, we obtain
\begin{equation}\label{e29}
\tilde{v}_{tt}+\sigma \tilde{v}_t-\gamma \tilde{v}+v^q_B=0.
\end{equation}
For problem (\ref{e29}), the following lemma is obtained in \cite{H}.

\begin{lem}\label{lemm9}
Let $v(t,\theta)$ be a positive axial symmetric solution of (\ref{e28})
satisfying $v_\theta>0$. Then $\tilde{v}(t)$ is bounded on $\mathbb{R}_+$.
\end{lem}

In order to prove that $v(t,\theta)$ is bounded on $\mathbb{R}_+\times (0,\frac{\pi}{2})$, we need to establish an apriori estimate for the following problem
\begin{equation}\label{e30}
-\mbox{div}(x_{n}^{a}\nabla u)+c(x)x_{n}^{a}u=0\;\;\mbox{in}\;\;B^+_1,~~~~~
\frac{\partial u}{\partial \nu^a}=F(x')u \;\;\;\mbox{on}\;\;\partial D_1,
\end{equation}
where $B_1^+:=\{x\in\mathbb{R}^n_+: |x|<1\}$, $D_1:=\{x\in\partial\mathbb{R}^n_+: |x|<1\}$.
We will borrow the method used in Theorem 8.17 in \cite{GT} (see also in \cite{H}) to prove the following  apriori estimate.

\begin{lem}\label{lemm10}
Let $u(x)\in H^1_a(B_1^+)$ be a weak solution of (\ref{e30})
with $F\in L^p(D_1)$ for some $p>\frac{n-1}{1-a}$. Then there exists $C>0$ depending on $n,~a,~ \|F\|_{L^p(D_1)}$ and $\|c\|_{L^\infty(B^+_{1})}$ such that
$$
\|u\|_{L^\infty(B^+_{1/2})}\leq C\|u\|_{L_a^2(B^+_{1})}.
$$
\end{lem}
\emph{Proof.} Let $\chi(|x|)$ be a smooth cut-off function satisfying $\chi(|x|)=0$ for $|x|\geq1$ and set $\tilde{c}=\|c\|_{L^\infty(B^+_{1})}$.
We denote $u_+=\max\{u,0\},~u_-=\max\{-u,0\}$.
Multiplying the equation in (\ref{e30}) by $u_+^k\chi^2(k>1)$ and integrating over $B^+_{1}$, we obtain
\begin{eqnarray}
\nonumber
&&\frac{4k}{(k+1)^2}\int_{B^+_{1}}|\nabla u_+^{(k+1)/2}|^2\chi^2x_{n}^{a}dx\\
\label{e31}
&&\leq
\int_{B^+_{1}}u_+^k|\nabla u_+||\nabla \chi^2|x_{n}^{a}dx
+\tilde{c}\int_{B^+_{1}}u_+^{k+1} \chi^2x_{n}^{a}dx
+\int_{D_{1}}F(x')u_+^{k+1} \chi^2dx'.
\end{eqnarray}
Let $p'=\frac{p}{p-1}$ be the dual exponent of $p$. Since $p>\frac{n-1}{1-a}$, then $2<2p'<\frac{2(n-1)}{n+a-2}$. We denote $2_\ast:=\frac{2(n-1)}{n+a-2}$.
Then there exists a $\varsigma\in(0,1)$ such that $\frac{1}{2p'}=\frac{\varsigma}{2}+\frac{1-\varsigma}{2_\ast}$.
By the H\"{o}lder inequality and an interpolation inequality, we have
\begin{eqnarray}
\nonumber
&&\int_{D_{1}}F(x')u_+^{k+1} \chi^2dx'\leq \|F\|_{L^p(D_1)}\|u_+^{(k+1)/2}\chi\|^2_{L^{2p'}(D_1)}\\
\nonumber
&&\leq
\|F\|_{L^p(D_1)}\left[\varepsilon\|u_+^{(k+1)/2}\chi\|_{L^{2_\ast}(D_1)}+\varepsilon^{-\frac{1-\varsigma}{\varsigma}}\|u_+^{(k+1)/2}\chi\|_{L^{2}(D_1)}\right]^2.
\end{eqnarray}
Since $a\in(-1,0)$, the following trace inequalities hold
$$
\|w\|_{L^{2_\ast}(D_1)}\leq C\|w\|_{H_a^1(B^+_{1})},~~~\|w\|_{L^{2}(D_1)}\leq A(\varepsilon,\varsigma)\|\nabla w\|_{L_a^2(B^+_{1})}+\frac{C}{A(\varepsilon,\varsigma)}\| w\|_{L_a^2(B^+_{1})},
$$
where $A(\varepsilon,\varsigma)>0$ is chosen small enough as $\varepsilon^{-\frac{1-\varsigma}{\varsigma}}A(\varepsilon,\varsigma)<\varepsilon$.
Hence we obtain
\begin{equation}\label{e32}
\int_{D_{1}}F(x')u_+^{k+1} \chi^2dx'\leq \|F\|_{L^p(D_1)}\left[4\varepsilon^2\|u_+^{(k+1)/2}\chi\|_{H_a^1(B^+_{1})}+C(\varepsilon)\|u_+^{(k+1)/2}\chi\|_{L_a^2(B^+_{1})}\right].
\end{equation}
From (\ref{e31}) and (\ref{e32}), choosing $\varepsilon^2=(\|F\|_{L^p(D_1)})^{-1}\frac{k}{2(k+1)^2}$, we obtain
$$
\int_{B^+_{1}}|\nabla (u_+^{(k+1)/2}\chi)|^2x_{n}^{a}dx\leq C(n,a,\|F\|_{L^p(D_1)},\tilde{c})(k+1)^2\int_{B^+_{1}} u_+^{k+1}(\chi+|\nabla \chi|)^2x_{n}^{a}dx.
$$
From \cite{Haj}, we know that there exists some $\upsilon>2$ such that for $2^\ast:=\frac{2\upsilon}{\upsilon-2}$, one has
\begin{equation}\label{e33}
\|u_+^{(k+1)/2}\chi\|_{L_a^{2^\ast}(B^+_{1})}\leq C(n,a, \|F\|_{L^p(D_1)},\tilde{c})(k+1)\|u_+^{(k+1)/2}(\chi+|\nabla \chi|)\|_{L_a^{2}(B^+_{1})}.
\end{equation}
It is now desirable to specify the cut-off function $\chi$ more precisely. Let $r_1,~r_2$ be such that $\frac{1}{2}\leq r_1\leq r_2\leq \frac{3}{4}$ and set
$\chi\equiv1$ in $B^+_{r_1}$, $\chi\equiv0$ in $B^+_{1}\backslash B^+_{r_2}$ with $|\nabla \chi|\leq \frac{2}{r_2-r_1}$.
Writing $\eta:=\frac{\upsilon}{\upsilon-2}$, we then have from (\ref{e33})
\begin{equation}\label{e34}
\|u_+^{(k+1)/2}\|_{L_a^{2\eta}(B^+_{r_1})}\leq C(k+1)\|u_+^{(k+1)/2}\|_{L_a^{2}(B^+_{r_2})}.
\end{equation}
Let us  introduce
$$
\Phi(\delta,r):=\left(\int_{B^+_{r}}| u_+|^\delta x_{n}^{a}dx\right)^{\frac{1}{\delta}}.
$$
Due to the fact that $\int_{B^+_{r}} x_{n}^{a}dx<\infty$, it is verified that
$$
\Phi(+\infty,r)=\lim_{\delta\rightarrow +\infty}\Phi(\delta,r)=\sup_{B^+_{r}} u_+,~~~~~
\Phi(-\infty,r)=\lim_{\delta\rightarrow -\infty}\Phi(\delta,r)=\inf_{B^+_{r}} u_+.
$$
From (\ref{e34}) we have
\begin{equation}\label{e35}
\Phi(\eta(k+1),r_1)\leq C(k+1)^{\frac{2}{(k+1)}}\Phi((k+1),r_2).
\end{equation}
This inequality can now be iterated to yield the desired estimates.
Hence, taking $\delta>1$, we set $k+1=\eta^i \delta$ and $r_i=2^{-1}+2^{-1-i},~i=0,1,\ldots$, so that, by (\ref{e35}) we have
$$
\Phi(\eta^i\delta,\frac{1}{2})\leq (C\eta)^{2\sum_i i\eta^{-i}}\Phi(\delta,1)=\tilde{C}\Phi(\delta,1),
$$
where $\tilde{C}$ depending on $\|F\|_{L^p(D_1)},\tilde{c}, a, n, \delta$.
Consequently, letting $i$ tend to infinity, we have
$$
\|u_+\|_{L^\infty(B^+_{1/2})}\leq C\|u_+\|_{L_a^\delta(B^+_{1})}.
$$
Let us now choose $\delta=2$, then 
$$
\|u_+\|_{L^\infty(B^+_{1/2})}\leq C\|u_+\|_{L_a^2(B^+_{1})}.
$$
By the same way, we can obtain the estimate for $u_-$.

 We complete the proof of this lemma.

\begin{mrem}\label{mrem2}
Similar estimate as that of Lemma \ref{lemm10} had been obtained by Fabes-Kenig-Serapioni \cite{FKS} (Corollary 2.3.4).
\end{mrem}

\begin{lem}\label{lemm11}
Let $v(t,\theta)$ be a positive axial symmetric solution of (\ref{e28})
satisfying $v_\theta>0$. Then $v(t,\theta)$ is bounded on $\mathbb{R}_+\times (0,\frac{\pi}{2})$.
\end{lem}
\emph{Proof.}
We apply a test function method to prove a boundedness of $v(t,\theta)$.
Let $\chi(t)$ be a cut-off function with a compact support in $\mathbb{R}_+$.
Multiplying (\ref{e28}) by $v^{-\frac{1}{2}}(t,\theta)\chi^2\sin^{n-2}\theta\cos^a\theta$ and integrating on $\mathbb{R}_+\times (0,\frac{\pi}{2})$,
we obtain
\begin{equation}\label{e36}
\int_{\mathbb{R}}\left(\int_0^{\frac{\pi}{2}}v^{-\frac{3}{2}}(v^2_t+v^2_\theta)d\mu+v_B^{q-1/2}\right)\chi^2dt
\leq C
\int_{\mathbb{R}}\int_0^{\frac{\pi}{2}}v^{\frac{1}{2}}(\chi^2+\chi^2_t)d\mu dt.
\end{equation}
We need to use the following trace inequality for a two dimensional domain $(a_1,a_2)\times (b_1,b_2)$
$$
\int_{a_1}^{a_2}|\phi(\xi_1,b_2)|^\alpha d\xi_1\leq C_\alpha\left(\int_{a_1}^{a_2}\int_{b_1}^{b_2}[\phi^2+\phi^2_{\xi_1}+\phi^2_{\xi_2}]d\xi_2d\xi_1\right)^{\frac{\alpha}{2}},
$$
where $\alpha\geq1$ and $C_\alpha>0$ is a constant depending on $\alpha$, $a_2-a_1$ and $b_2-b_1$.
Applying this inequality with $\phi(\xi_1, \xi_2)=v^{\frac{1}{4}}(\xi_1, \xi_2)$, $a_1=\tau-1,~a_2=\tau+1,~b_1=\frac{\pi}{4},~b_1=\frac{\pi}{2}$,
then we obtain
$$
\int_{\tau-1}^{\tau+1}|v^{\frac{1}{4}}_B|^\alpha dt\leq C_\alpha\left(\int_{\tau-1}^{\tau+1}\int_{\frac{\pi}{4}}^{\frac{\pi}{2}}[v^{\frac{1}{2}}+v^{-\frac{3}{2}}(v^2_t+v^2_\theta)]d\theta dt\right)^{\frac{\alpha}{2}}.
$$
We take a cut-off function $\chi(t)$ such that $\chi(t)=1$ if $t\in [\tau-1,\tau+1]$ and $\chi(t)=0$ if $t\in \mathbb{R}\backslash (\tau-2,\tau+2)$.
It is clear that $\sin\theta\geq\frac{\sqrt{2}}{2}$ for $\theta\in[\frac{\pi}{4},\frac{\pi}{2}]$, and $\cos^a\theta>1$, since $a\in(-1,0)$.

Hence from (\ref{e36}), for $\alpha\geq1$, we have
$$
\int_{\tau-1}^{\tau+1}|v^{\frac{1}{4}}_B|^\alpha dt\leq C_\alpha\left(\int_{\tau-2}^{\tau+2}\int_{0}^{\frac{\pi}{2}}v^{\frac{1}{2}}d\mu dt\right)^{\frac{\alpha}{2}}\leq C_\alpha\left(\int_{\tau-2}^{\tau+2}\tilde{v} dt\right)^{\frac{\alpha}{4}},
$$
where in the last inequality we have used H\"{o}lder inequality.
Therefore, from Lemma \ref{lemm9}, for $\alpha\geq1$ there exists $\bar{C}_\alpha>0$ independent of $\tau>0$ such that
$$
\int_{\tau-1}^{\tau+1}|v^{\frac{1}{4}}_B|^\alpha dt\leq \bar{C}_\alpha.
$$
Then from Lemma \ref{lemm10}, we deduce that $v(t,\theta)$ is bounded on $\mathbb{R}_+\times (0,\frac{\pi}{2})$.

The proof of this lemma is complete.

In the rest of this section, we will prove  Theorem \ref{thm1} (iii): the limit behavior of positive axial symmetric solutions.

Recall that $v(t,\theta)$ satisfies (\ref{e28}). Note that $\sigma,~\gamma>0$ if $q>\frac{n-a}{n+a-2}$.
We define the following energy function $E(t)$ associated with (\ref{e28}) by
$$
E(t):=\frac{1}{2}\|\partial_tv\|^2_{2,a}-\frac{\gamma}{2}\|v\|^2_{2,a}-\frac{1}{2}\|\partial_\theta v\|^2_{2,a}+\frac{1}{q+1}(v_B(t))^{q+1}.
$$
Simple calculation shows that
$$
\partial_tE(t)=-\sigma\|\partial_tv\|^2_{2,a}.
$$

\begin{lem}\label{lemm12}
Let $q>\frac{n-a}{n+a-2}$ and $v(t,\theta)$ be a positive bounded solution of (\ref{e28})
with $\lim_{t\rightarrow-\infty}E(t)=0$. Then one has
$$\lim_{t\rightarrow+\infty}v(t,\theta)=V(\theta)~~\mbox{in}~C[0,\pi/2].$$
\end{lem}
\emph{Proof.} From the boundary condition in (\ref{e28}), $a\in(-1,0)$ and the result that $v(t,\theta)$ is bounded on $\mathbb{R}\times(0,\pi/2)$,
 we have $ v_\theta(t,\frac{\pi}{2})=0$. The elliptic regularity theory assures a
boundedness of $\partial_tv(t,\theta)$ and $\partial_\theta v(t,\theta)$. Hence we deduce that $E(t)$ is bounded and
\begin{equation}\label{e37}
0<\sigma\int_{-\infty}^\infty\int_0^{\frac{\pi}{2}}v^2_td\mu dt=\lim_{t\rightarrow-\infty}E(t)-\lim_{t\rightarrow\infty}E(t)=-\lim_{t\rightarrow\infty}E(t)<\infty.
\end{equation}
Let $\{t_i\}_{i\in N}$ be any sequence satisfying $\lim_{i\rightarrow\infty}t_i=+\infty$ and set $v_i(t,\theta)=v(t+t_i, \theta)$.
Then there exists a subsequence of $\{t_i\}_{i\in N}$, which is still denoted by the same symbol such that $v_i(t,\theta)$
converges to some function $v_\infty(t,\theta)$ in $C([-1,1]\times[0,\pi/2])$.
By (\ref{e37}) we find that $\partial_tv_\infty(t,\theta)\equiv0$. Then we denote $v_\infty(t,\theta)=:v_\infty(\theta)$ and $v_\infty(\theta)$
satisfies
$$
-\gamma v+\Delta_{S,a}v=0~~~\mbox{in}~(0,\pi/2),~~~~\lim_{\theta\rightarrow \frac{\pi}{2}}v_\theta\cos^a \theta =v_B^{q}.
$$
By $\lim_{t\rightarrow-\infty}E(t)=0$ and (\ref{e37}), we have $\lim_{t\rightarrow\infty}E(t)<0$, which yields that the limiting function $v_\infty(\theta)$
is not a trivial one.
Hence by Proposition \ref{prop2}, we deduce that $v_\infty(\theta)\equiv V(\theta)$.
Therefore we obtain the result that $\lim_{i\rightarrow\infty}v(t_i,\theta)=V(\theta)$ in $C[0,\pi/2]$ for any sequence $\{t_i\}_{i\in N}$ converging to $+\infty$. Hence
$$\lim_{t\rightarrow\infty}v(t,\theta)=V(\theta)~~~~~~\mbox{in}~C[0,\pi/2].$$
The lemma is proved.

\emph{Proof of  Theorem \ref{thm1} (iii)}
Recalling that $\psi_\infty(x)=V(\theta)r^{\frac{a-1}{q-1}}$ and $u_\beta(x)=\beta u_1(\beta^{\frac{q-1}{1-a}}x)$,
we have
$$
|u_\beta(x)-\psi_\infty(x)|=\left|\beta u_1(\beta^{\frac{q-1}{1-a}}x)-V(\theta)|x|^{\frac{a-1}{q-1}}\right|
=\left|v_1(\beta^{\frac{q-1}{1-a}}x)-V(\theta)\right||x|^{\frac{a-1}{q-1}}.
$$
Note that $\lim_{t\rightarrow-\infty}v(t,\theta)=0$, since $u$ is bounded.
Then the condition $\lim_{t\rightarrow-\infty}E(t)=0$ in Lemma \ref{lemm12} is true.
Hence by Lemma \ref{lemm12}, we have
$$
\lim_{\beta\rightarrow+\infty}|u_\beta(x)-\psi_\infty(x)|=0,~~~x\neq0.
$$
We complete the proof of Theorem \ref{thm1} (iii).

\section{Asymptotic expansion}

In this section, we will prove Theorem \ref{thm2}.
We need to study the asymptotic behavior of
$$
w(t,\theta)=V(\theta)-v(t,\theta).
$$
Note that $w(t,\theta)$ satisfies
\begin{equation}\label{e38}
\left\{\begin{array}{l}w_{tt}+\sigma w_t-\gamma w+\Delta_{S,a}w=0\;\;\;\;\;~~~~~~~~~~\mbox{in}\;\;\mathbb{R}\times (0,\frac{\pi}{2}),\\
\lim_{\theta\rightarrow \frac{\pi}{2}}w_\theta\cos^a \theta =qv^{q-1}_Bw_B+g(w_B) \;\;\;\;\;\mbox{on}\;\; \mathbb{R},\\
\end{array}
\right.
\end{equation}
where $g(w)$ is given by
$$
g(w)=v^{q}_B-(v_B-w)^{q}-qv^{q-1}_Bw.
$$
Since $L^2_{\mbox{sym}}(S^{n-1}_+)$ is spanned by the eigenfunctions $\{e_i(\theta)\}_{i\in N}$ of (\ref{e23}), then $w(t,\theta)$
can be expanded as
$$
w(t,\theta)=\sum_{i=1}^\infty z_i(t)e_i(\theta),~~t\in\mathbb{R}.
$$
Multiplying the equation in (\ref{e23}) by $e_i(\theta)\sin^{n-2}\theta\cos^a\theta$ and integrating with respect to $\theta$
on $(0,\frac{\pi}{2})$, we obtain
\begin{equation}\label{e39}
z_i''+\sigma z_i'-(\gamma+\lambda_i)z_i=f_i(t),
\end{equation}
where $f_i(t):=-g(w_B(t))e_{iB}$.

We first consider the case $i\geq2$.
Recall that Lemma \ref{lemm8} shows that $\lambda_i>0$ for $i\geq2$.
Hence the corresponding quadratic equation
\begin{equation}\label{e40}
\rho^2+\sigma \rho-(\gamma+\lambda_i)=0
\end{equation}
to (\ref{e39}) admits two real roots
\begin{equation}\label{e41}
\rho_i^{\pm}=\frac{-\sigma\pm\sqrt{\sigma^2+4(\gamma+\lambda_i)}}{2}.
\end{equation}
Note that $\rho_i^{-}<0<\rho_i^{+}.$
From \cite{H} we know that for $i\geq2$
\begin{eqnarray}
\nonumber z_i(t)&=& z_i(0)e^{\rho_i^{-}t}-\frac{e^{\rho_i^{-}t}}{\sqrt{\sigma^2+4(\gamma+\lambda_i)}}\int_0^t(e^{-\rho_i^{-}s}-e^{-\rho_i^{+}s})f_i(s)ds\\
\label{e42} &&-\frac{e^{\rho_i^{+}t}-e^{\rho_i^{-}t}}{\sqrt{\sigma^2+4(\gamma+\lambda_i)}}\int_t^\infty e^{-\rho_i^{+}s}f_i(s)ds.
\end{eqnarray}

For the case $i=1$, by Lemma \ref{lemm7}, (\ref{e39}) admits two real roots
 $\rho_1^{-}<\rho_1^{+}<0$ if $q$ is JL-supercritical, just one root $\rho_1=-\frac{\sigma}{2}$ if $q$ is JL-critical and admits no real roots if
$q$ is JL-subcritical. Hence we obtain for the JL-supercritical case
\begin{eqnarray}
\nonumber z_1(t)&=& z_1(0)e^{\rho_1^{-}t}-\frac{z'_1(0)-\rho_1^{-}z_1(0)}{\sqrt{\sigma^2+4(\gamma+\lambda_1)}}(e^{\rho_1^{+}t}-e^{\rho_1^{-}t})\\
\label{e43} &&+\int_0^t \frac{1-e^{(2\rho_1^{-}+\sigma)(t-s)}}{\sqrt{\sigma^2+4(\gamma+\lambda_1)}}e^{\rho_1^{+}(t-s)}f_1(s)ds,
\end{eqnarray}
for the JL-critical case
\begin{equation}\label{e44}
z_1(t)= z_1(0)e^{\rho_1t}+(z'_1(0)-\rho_1z_1(0))te^{\rho_1t}+\int_0^t (t-s)e^{\rho_1(t-s)}f_1(s)ds
\end{equation}
and for the JL-subcritical case
\begin{eqnarray}
\nonumber z_1(t)&=&\frac{1}{K} \left(\frac{\sigma}{2}z_1(0)+z'_1(0)\right)(\sin Kt)e^{-\frac{\sigma}{2}t}+z_1(0)(\cos Kt)e^{-\frac{\sigma}{2}t}
\\
\label{e45} &&+\frac{1}{K}\int_0^t [(\sin Kt)(\cos Ks)-(\sin Ks)(\cos Kt)]e^{-\frac{\sigma(t-s)}{2}}f_1(s)ds,
\end{eqnarray}
where $K$, given in  Theorem \ref{thm2}, is the imaginary part of a root of $\rho^2+\sigma \rho-(\gamma+\lambda_1)=0$.

To prove  Theorem \ref{thm2}, we need to establish the following  proposition.

\begin{pro} \label{prop3}
(i) If $q$ is JL-supercritical, then there exists $\xi_1>0$ and $\varepsilon,~C>0$ such that
$$
\|w(t,\theta)-\xi_1e^{\rho_1^{+}t}e_1\|_\infty\leq Ce^{(\rho_1^{+}-\varepsilon)t},~~~t>0;
$$

(ii) If $q$ is JL-critical, then there exists $\xi_1>0,~\xi_2\in \mathbb{R}$ and $\varepsilon,~C>0$ such that
$$
\|w(t,\theta)-(\xi_1t+\xi_2)e^{-\sigma t/2}e_1\|_\infty\leq Ce^{-(\sigma+\varepsilon)t/2},~~~t>0;
$$

(iii) If $q$ is JL-subcritical, then one of the following two expansions holds.

(iii-1) there exist $(\xi_1,\xi_2)\neq0$ and $\varepsilon,~C>0$ such that
$$
\|w(t,\theta)-(\xi_1(\sin Kt)+\xi_2(\cos Kt))e^{-\sigma t/2}e_1\|_\infty\leq Ce^{-(\sigma+\varepsilon)t/2},~~~t>0;
$$

(iii-2) there exist $\xi\neq0$ and $\varepsilon,~C>0$ such that
$$
\|w(t,\theta)-\xi e^{\rho_2^-t}e_2\|_\infty\leq Ce^{(\rho_2^{-}-\varepsilon)t},~~~t>0.
$$
\end{pro}

As a consequence of Proposition \ref{prop3}, we immediately obtain Theorem \ref{thm2}.

\emph{Proof of Theorem \ref{thm2}.} We only show the proof of JL-supercritical case in Theorem \ref{thm2}, since the proof of the rest part  is similar.
Let  $q$ be JL-supercritical and set $\rho_1=|\rho_1^{+}|$.
By Proposition \ref{prop3} we obtain
$$
\|v(t,\theta)-(V-\xi_1e^{-\rho_1t}e_1)\|_\infty\leq Ce^{-(\rho_1+\varepsilon)t},~~~t>0.
$$
Going back to the original variable, we obtain
$$
u(r,\theta)=V(\theta)r^{-m_q}-\xi_1r^{-m_q-\rho_1}e_1(\theta)+\mbox{O}(r^{-m_q-\rho_1-\varepsilon}),~~~r>>1.
$$
Plugging the explicit expression of $\sigma$ and $\gamma$ into (\ref{e41}), we obtain the explicit expression of $\rho_1$,
which coincides with the expression in Theorem \ref{thm2}.

Now the remaining  task for us is to prove Proposition \ref{prop3}.

First we show that $w(t,\theta)$ decays exponentially as $t\rightarrow+\infty$.

\begin{lem}\label{lemm13}
There exists $\varepsilon,~C>0$ such that
$$
\|w(t,\cdot)\|_{2,a}\leq Ce^{-\varepsilon t},~~~t>0.
$$
\end{lem}
\emph{Proof.}
From (\ref{e42}), we obtain for $i\geq2$
\begin{eqnarray}
\nonumber |z_i(t)|&\leq & |z_i(0)|e^{\rho_i^{-}t}+\frac{e^{\rho_i^{-}t}}{\sqrt{\lambda_i}}\int_0^te^{-\rho_i^{-}s}|f_i(s)|ds
+\frac{e^{\rho_i^{+}t}}{\sqrt{\lambda_i}}\int_t^\infty e^{-\rho_i^{+}s}|f_i(s)|ds\\
\nonumber  &\leq&|z_i(0)|e^{\rho_2^{-}t}+\frac{e^{\rho_i^{-}t/2}}{\sqrt{\lambda_i|\rho_i^{-}|}}\left(\int_0^te^{-\rho_i^{-}s}|f_i(s)|^2ds\right)^{1/2}\\
\nonumber  &&+\frac{e^{\rho_i^{+}t/2}}{\sqrt{\lambda_i|\rho_i^{+}|}}\left(\int_t^\infty e^{-\rho_i^{+}s}|f_i(s)|^2ds\right)^{1/2}.
\end{eqnarray}
By a trace inequality
$$
|e_{iB}|\leq C\|e_i\|^{1/2}_{2,a}\|\partial_\theta e_i\|^{1/2}_{2,a},
$$
we have
\begin{eqnarray}
\nonumber \|\partial_\theta e_i\|^2_{2,a}&= & 2\|\partial_\theta e_i\|^2_{2,a}-\|\partial_\theta e_i\|^2_{2,a}
=2(qV_B^{q-1}e_{iB}^2+\lambda_i\|e_i\|^2_{2,a})-\|\partial_\theta e_i\|^2_{2,a}\\
\nonumber  &\leq&2\lambda_i\|e_i\|^2_{2,a}+CqV_B^{q-1}\| e_i\|_{2,a}\|\partial_\theta e_i\|_{2,a}-\|\partial_\theta e_i\|^2_{2,a}\\
\nonumber  &\leq&2\lambda_i\|e_i\|^2_{2,a}+C\| e_i\|^2_{2,a}\leq C\lambda_i\|e_i\|^2_{2,a}.
\end{eqnarray}
Therefore, since $\|e_i\|_{2,a}=1$, we deduce that $\|\partial_\theta e_i\|^2_{2,a}\leq C\lambda_i$. Then we have
$$
|e_{iB}|\leq C\|e_i\|^{1/2}_{2,a}\|\partial_\theta e_i\|^{1/2}_{2,a}\leq C\lambda_i^{\frac{1}{4}}.
$$
Recalling that $f_i(t)=-g(w_B(t))e_{iB}$, we obtain
$$
|f_i(s)|\leq C\lambda_i^{\frac{1}{4}}w^2_B(s)\leq C\lambda_i^{\frac{1}{4}}\|w(s,\cdot)\|_{2,a}\|\partial_\theta w(s,\cdot)\|_{2,a}.
$$
Since $|\rho_i^{-}|,~|\rho_i^{+}|\geq C\sqrt{\lambda_i}$ for some $C>0$, we have for $i\geq2$
\begin{eqnarray}
\nonumber |z_i(t)|^2&\leq & C|z_i(0)|^2e^{2\rho_2^{-}t}+\frac{C}{\lambda_i}(\int_0^te^{\rho_i^{-}(t-s)}\varphi(s)\|w(s,\cdot)\|^2_{2,a}ds
\\
\nonumber  & &
+\int_t^\infty e^{-\rho_i^{+}(s-t)}\varphi(s)\|w(s,\cdot)\|^2_{2,a}ds)\\
\nonumber  &\leq&C|z_i(0)|^2e^{2\rho_2^{-}t}+\frac{C}{\lambda_i}(\int_0^te^{\rho_2^{-}(t-s)}\varphi(s)\|w(s,\cdot)\|^2_{2,a}ds
\\
\nonumber  & &
+\int_t^\infty e^{-\rho_2^{+}(s-t)}\varphi(s)\|w(s,\cdot)\|^2_{2,a}ds),
\end{eqnarray}
where $\varphi(s):=\|\partial_\theta w(s,\cdot)\|^2_{2,a}$. Applying Strum's comparison theorem, we can obtain the following estimate of the eigenvalues of problem (\ref{e23})
$$
C_1i^2\leq \lambda_i\leq C_2i^2,~~~\forall ~i>1,
$$
where $C_1,~C_2$ are positive constants independence of $i$.
Hence we have
\begin{eqnarray}
\nonumber \sum_{i=2}^\infty|z_i(t)|^2&\leq & C\|w(0,\cdot)\|^2_{2,a}e^{2\rho_2^{-}t}+C(\int_0^te^{\rho_2^{-}(t-s)}\varphi(s)\|w(s,\cdot)\|^2_{2,a}ds
\\
\nonumber  & &
+\int_t^\infty e^{-\rho_2^{+}(s-t)}\varphi(s)\|w(s,\cdot)\|^2_{2,a}ds).
\end{eqnarray}
Similarly, from (\ref{e43})-(\ref{e45}), we obtain
\begin{eqnarray}\nonumber
|z_1(t)|^2\leq\left\{\begin{array}{l}Ce^{2\rho_1^{+}t}+C\int_0^te^{\rho_1^{+}(t-s)}\varphi(s)\|w(s,\cdot)\|^2_{2,a}ds\;\;\;\;\;~~~~~~~~~~
\mbox{if}\;q~ \mbox{is JL-supercritical},\\
Ct^2e^{2\rho_1t}+C\int_0^t(t-s)^2e^{\rho_1(t-s)}\varphi(s)\|w(s,\cdot)\|^2_{2,a}ds\;\;\;\;\;
\mbox{if}\;q~ \mbox{is JL-critical},\\Ce^{-\sigma t}+C\int_0^te^{-\sigma(t-s)}\varphi(s)\|w(s,\cdot)\|^2_{2,a}ds\;\;\;\;\;~~~~~~~~~~~
\mbox{if}\;q~ \mbox{is JL-subcritical}.\\
\end{array}
\right.
\end{eqnarray}
Hence there exists $\delta>0$ such that
$$
\|w(t,\cdot)\|^2_{2,a}\leq Ce^{-\delta t}+C\left(\int_0^te^{-\delta(t-s)}\varphi(s)\|w(s,\cdot)\|^2_{2,a}ds
+\int_t^\infty e^{-\delta(s-t)}\varphi(s)\|w(s,\cdot)\|^2_{2,a}ds\right)
$$
By Lemma  \ref{lemm14} stated below, we deduce the desired conclusion.

\begin{lem}\label{lemm14}(\cite{H})
Assume $\eta(t)$ and $\varphi(t)$ are positive continuous functions defined on $\overline{\mathbb{R}}_{+}$
converging to zero as $t\rightarrow+\infty$. Moreover $\eta(t)$ satisfies
$$
\eta(t)\leq Ce^{-\delta t}+C\left(\int_0^te^{-\delta(t-s)}\varphi(s)\eta(s)ds
+\int_t^\infty e^{-\delta(s-t)}\varphi(s)\eta(s)ds\right),~~~t>0
$$
for some $\delta>0$. Then there exists $\varepsilon>0$ such that
$
\eta(t)\leq Ce^{-\varepsilon t}.
$
\end{lem}

Next we will show more precise decay rates of $\|w(t,\cdot)\|_{2,a}$.

\begin{lem}\label{lemm15}
There exists $C>0$ such that for $t>0$
\begin{eqnarray}\nonumber
\|w(t,\cdot)\|_{2,a}\leq\left\{\begin{array}{l} Ce^{\rho_1^{+}t}\;\;\;\;\;~~~~~~~~~~
\mbox{if}\;q~ \mbox{is JL-supercritical},\\
C(1+t)e^{\rho_1t}\;\;\;\;\;~~
\mbox{if}\;q~ \mbox{is JL-critical},\\Ce^{-\sigma t/2}\;\;\;\;\;~~~~~~~
\mbox{if}\;q~ \mbox{is JL-subcritical}.\\
\end{array}
\right.
\end{eqnarray}
\end{lem}
\emph{Proof.} Repeating the argument given in the proof of Lemma \ref{lemm13}, we have
\begin{eqnarray}
\nonumber \sum_{i=k}^\infty|z_i(t)|^2&\leq & Ce^{2\rho_k^{-}t}+C(\int_0^te^{\rho_k^{-}(t-s)}\varphi(s)\|w(s,\cdot)\|^2_{2,a}ds
\\
\nonumber  & &
+\int_t^\infty e^{-\rho_k^{+}(s-t)}\varphi(s)\|w(s,\cdot)\|^2_{2,a}ds).
\end{eqnarray}
Since $w(t,\cdot)$ is uniformly bounded in $\mathcal{H}^2_a(S^{n-1}_+)$,
by an interpolation inequality, it holds that
$$
\|\partial_\theta w(t,\cdot)\|^2_{2,a}\leq C\| w(t,\cdot)\|_{2,a}\| w(t,\cdot)\|_{\mathcal{H}^2_a(S^{n-1}_+)}\leq C\| w(t,\cdot)\|_{2,a}.
$$
Recalling $\varphi(s)=\|\partial_\theta w(s,\cdot)\|^2_{2,a}$, by Lemma \ref{lemm13}, we have for $k\geq2$
$$
\sum_{i=k}^\infty|z_i(t)|^2\leq  C[e^{2\rho_k^{-}t}+e^{\rho_k^{-}t}+e^{-3\varepsilon t}].
$$
From (\ref{e42})-(\ref{e45}), we have that for $2\leq i\leq k-1$
\begin{eqnarray}
\nonumber |z_i(t)|&\leq & Ce^{\rho_i^{-}t}+Ce^{\rho_i^{-}t}\int_0^te^{-\rho_i^{-}s}e^{-3\varepsilon s/2}ds
+Ce^{\rho_i^{+}t}\int_t^\infty e^{-\rho_i^{+}s}e^{-3\varepsilon s/2}ds\\
\nonumber  &\leq&C[e^{\rho_i^{-}t}+e^{-3\varepsilon t/2}]\leq C[e^{\rho_2^{-}t}+e^{-3\varepsilon t/2}]
\end{eqnarray}
and
\begin{eqnarray}\label{e46}
|z_1(t)|\leq\left\{\begin{array}{l}C[e^{\rho_1^{+}t}+e^{-3\varepsilon t/2}]\;\;\;\;\;~~~~~~~~~~
\mbox{if}\;q ~\mbox{is JL-supercritical},\\
C[(1+t)e^{\rho_1t}+e^{-3\varepsilon t/2}]~\;\;\;\;\;~
\mbox{if}\;q~ \mbox{is JL-critical},\\C[e^{-\sigma t/2}+e^{-3\varepsilon t/2}]\;\;\;\;\;~~~~~~~
\mbox{if}\;q~ \mbox{is JL-subcritical}.\\
\end{array}
\right.
\end{eqnarray}
We choose $k\geq2$ large enough such that $|\rho_k^-|\geq2|\rho_2^-|$, then we obtain
\begin{equation}\label{e47}
\sum_{i=2}^\infty|z_i(t)|^2\leq  C[e^{2\rho_2^{-}t}+e^{-3\varepsilon t}].
\end{equation}

We first consider the JL-supercritical case.
If $\varepsilon\geq \frac{2|\rho_1^+|}{3}$, from (\ref{e46})-(\ref{e47}) we know that  this lemma is true. If $\varepsilon< \frac{2|\rho_1^+|}{3}$,
from (\ref{e46})-(\ref{e47}) we know that
$$
\| w(t,\cdot)\|_{2,a}\leq Ce^{-3\varepsilon t/2}.
$$
Repeating the above procedure, if $\frac{2|\rho_1^+|}{3}>\varepsilon\geq \frac{4|\rho_1^+|}{9}$, from (\ref{e46})-(\ref{e47}) (now the quantity $\varepsilon$
in these inequalities should be changed into the quantity $\frac{3\varepsilon}{2}$), we know that  this lemma is also true. If $\varepsilon< \frac{4|\rho_1^+|}{9}$,
 we know that
$$
\| w(t,\cdot)\|_{2,a}\leq Ce^{-9\varepsilon t/4}.
$$
After repeating such procedure  finite times, we can obtain
$$
\| w(t,\cdot)\|_{2,a}\leq Ce^{-\rho_1^+t}.
$$
Therefore we complete the proof of the JL-supercritical case.
The JL-critical case and JL-subcritical case can be proved similarly.

Furthermore we see that the decay rate of $\| w(t,\cdot)\|_{2,a}$ given in Lemma \ref{lemm15} is exact one.

\begin{lem}\label{lemm16}
There exists $\xi_1,~\xi_2\in \mathbb{R}$ and $\varepsilon,~C>0$ such that for $t>0$
\begin{eqnarray}\nonumber
&&\|e^{-\rho_1^+t}w(t,\cdot)-\xi_1e_1\|_{2,a}\leq Ce^{-\varepsilon t}\;\;\;\;\;~~~~~~~~~~~~~~~~~~~~~~~~~~
\mbox{if}\;q \mbox{is JL-supercritical},\\
\label{e48}
&&\|e^{\sigma t/2}w(t,\cdot)-(\xi_1t+\xi_2)e_1\|_{2,a}\leq Ce^{-\varepsilon t}\;\;\;\;\;~~~~~~~~~~~~~~~~~~
\mbox{if}\;q \mbox{is JL-critical},\\
\nonumber
&&\|e^{\sigma t/2}w(t,\cdot)-(\xi_1 \sin (Kt)+\xi_2\cos(Kt))e_1\|_{2,a}\leq Ce^{-\varepsilon t}\;\;\;\;
\mbox{if}\;q \mbox{is JL-subcritical}.
\end{eqnarray}
Moreover the norm $\| \cdot\|_{2,a}$ can be replaced by the norm $\| \cdot\|_{\infty}$.
\end{lem}
\emph{Proof.}
By Lemma \ref{lemm15} and (\ref{e47}), for the  JL-supercritical case, we have
$$
\sum_{i=2}^\infty|z_i(t)|^2\leq  C[e^{2\rho_2^{-}t}+e^{3\rho_1^+t}].
$$
Hence it is sufficient to estimate $z_1(t)$.
Similarly, we also only need to estimate $z_1(t)$ in the JL-critical and JL-subcritical case.

We first consider the JL-supercritical case.
We set $z_1(t)=e^{\rho_1^+t}y_1(t)$.
From (\ref{e39}),  $y_1(t)$ solves
\begin{equation}\label{e49}
y''_1+(2\rho_1^++\sigma)y'_1=e^{-\rho_1^+t}f_1.
\end{equation}
Then since $e^{(2\rho_1^++\sigma)t}[y''_1+(2\rho_1^++\sigma)y'_1]=(e^{(2\rho_1^++\sigma)t}y'_1)'=e^{(\rho_1^++\sigma)t}f_1$, we have
$$
y'_1(t)=y'_1(0)e^{-(2\rho_1^++\sigma)t}+e^{-(2\rho_1^++\sigma)t}\int_0^te^{(\rho_1^++\sigma)s}f_1(s)ds.
$$
Hence by Lemma \ref{lemm15}, there exists $\varepsilon>0$ such that $|y'_1(t)|\leq Ce^{-\varepsilon t}$, which yields that
$\lim_{t\rightarrow +\infty}y_1(t)$ exists. We denote $\xi_1:=\lim_{t\rightarrow +\infty}y_1(t)$.
So we have
$$
|y_1(t)-\xi_1|\leq \int_t^\infty|y'_1(s)|ds\leq Ce^{-\varepsilon t},
$$
which shows (\ref{e48}) for the JL-supercritical case.

For the JL-critical case, we set $z_1(t)=te^{\rho_1t}y_1(t)$. From (\ref{e39}),  $y_1(t)$ satisfies
\begin{equation}\label{e50}
ty''_1+2y'_1=e^{-\rho_1t}f_1.
\end{equation}
This and the equation $t^2y''_1+2ty'_1=(t^2y'_1)'$ give that
\begin{eqnarray}\nonumber
y'_1(t)&=&\frac{1}{t^2}\int_0^tse^{-\rho_1s}f_1(s)ds\\
\label{e51}
&=&\frac{1}{t^2}\int_0^\infty se^{-\rho_1s}f_1(s)ds-\frac{1}{t^2}\int_t^\infty se^{-\rho_1s}f_1(s)ds,
\end{eqnarray}
which gives that $y'_1\in L^1([1,\infty))$. So $\lim_{t\rightarrow +\infty}y_1(t)$ exists and we denote $\xi_1:=\lim_{t\rightarrow +\infty}y_1(t)$.
We set $\xi_2:=\int_0^\infty se^{-\rho_1s}f_1(s)ds$.
Then, from (\ref{e51}), we have
$$
\left|y_1(t)-\xi_1+\frac{\xi_2}{t}\right|\leq  Ce^{-\varepsilon t},~~t>0
$$
for some $\varepsilon>0$, which shows (\ref{e48}) for the JL-critical case.

Finally we consider the  JL-subcritical case.
From (\ref{e45}), we have
\begin{eqnarray}
\nonumber z_1(t)&=&\frac{1}{K} \left(\frac{\sigma}{2}z_1(0)+z'_1(0)\right)(\sin Kt)e^{-\frac{\sigma}{2}t}+z_1(0)(\cos Kt)e^{-\frac{\sigma}{2}t}
\\
\nonumber  &&+\frac{e^{-\sigma t/2}}{K}\left((\sin Kt)\int_0^t (\cos Ks)e^{\sigma s/2}f_1(s)ds-(\cos Kt)\int_0^t (\sin Ks)e^{\sigma s/2}f_1(s)ds\right).
\end{eqnarray}
We set $\alpha_1:=\int_0^t (\cos Ks)e^{\sigma s/2}f_1(s)ds$ and $\alpha_2:=\int_0^t (\sin Ks)e^{\sigma s/2}f_1(s)ds$.
From Lemma \ref{lemm15}, we deduce that $|\alpha_1|,~|\alpha_2|<+\infty$.
Therefore we obtain
\begin{eqnarray}
\nonumber &&\left|z_1(t)-\frac{1}{K} \left(\frac{\sigma}{2}z_1(0)+z'_1(0)+\alpha_1\right)(\sin Kt)e^{-\frac{\sigma}{2}t}
-\left(z_1(0)-\frac{\alpha_2}{K}\right)(\cos Kt)e^{-\frac{\sigma}{2}t}\right|\\
\label{e52}
&&\leq\frac{1}{K}e^{-\frac{\sigma}{2}t}\int_t^\infty e^{\frac{\sigma}{2}s}|f_1(s)|ds,
\end{eqnarray}
which shows (\ref{e48}) for the JL-subcritical case.

Now we consider the $L^\infty$-estimate.
For the JL-supercritical case,
we set $Y(t,\theta):=e^{-\rho_1^+t}w(t,\theta)-\xi_1e_1(\theta)$.
Then we have
\begin{equation}\label{e53}
\left\{\begin{array}{l}Y_{tt}+(2\rho_1^++\sigma) Y_t+\lambda_1 Y+\Delta_{S,a}Y=0\;\;\;\;\;~~~\mbox{in}\;\;\mathbb{R}\times (0,\frac{\pi}{2}),\\
\lim_{\theta\rightarrow \frac{\pi}{2}}Y_\theta\cos^a \theta =qv^{q-1}_BY_B+e^{-\rho_1^+t}g(w_B) \;\;\;\;\mbox{on}\;\; \mathbb{R}.\\
\end{array}
\right.
\end{equation}
By the similar argument as given in Lemma \ref{lemm10}, we know that there exists $C>0$ such that
$$
\sup_{t-1\leq\tau\leq t+1}\|Y(\tau,\cdot)\|_{\infty}\leq C\left(\int_{t-2}^{t+2}\|Y(\tau,\cdot)\|^2_{2,a}d\tau\right)^{\frac{1}{2}}+C
\sup_{t-2\leq\tau\leq t+2}e^{-\rho_1^+\tau}|g(w_B(\tau))|.
$$
Since $|g(w_B(\tau))|\leq Cw^2_B(\tau)$, by a trace inequality and an interpolation, we have
$$
|g(w_B(t))|\leq Cw^2_B(t)\leq C\|\partial_\theta w(t,\cdot)\|_{2,a}\| w(t,\cdot)\|_{2,a}
\leq C\| w(t,\cdot)\|^{\frac{1}{2}}_{\mathcal{H}^2_a(S^{n-1}_+)}\| w(t,\cdot)\|^{\frac{3}{2}}_{2,a}.
$$
Therefore, since $\| w(t,\cdot)\|_{\mathcal{H}^2_a(S^{n-1}_+)}$ is bounded, we deduce that
$|g(w_B(t))|\leq Ce^{3\rho_1^{+} t/2}$, which yields the $L^\infty$-estimate
for the JL-supercritical case.

For the JL-critical and JL-subcritical case, the argument is similar, we omit them.

The proof of this lemma is complete.

\begin{lem}\label{lemm17}
Let $q$ be JL-supercritical, then $\xi_1>0$, where $\xi_1$ is the constant given in the previous lemma.
\end{lem}
\emph{Proof.} Recalling that $\lim_{t\rightarrow +\infty}y_1(t)=\xi_1$, we know that there exists a sequence $\{t_i\}_{i\in N}$ such that
$$
\lim_{i\rightarrow\infty}y_1(t_i)=\xi_1,~~~\lim_{i\rightarrow\infty}y'_1(t_i)=0.
$$
Integrating  (\ref{e49}) on $(t,t_i)$ and taking a limit $i\rightarrow\infty$, we obtain
$$
-y'_1(t)+(2\rho_1^++\sigma)(\xi_1-y_1(t))=\int_t^\infty e^{-\rho_1^+s}f_1(s)ds.
$$
Note that both $y_1(t)$ and $y'_1(t)$ converge to zero as $t\rightarrow-\infty$, we deduce that
$$
(2\rho_1^++\sigma)\xi_1=\int_{-\infty}^\infty e^{-\rho_1^+s}f_1(s)ds.
$$
By the convexity of the function $w^q$, we have $-g(w)=(V_B-w)^q-[V_B^q+qV_B^{q-1}(-w)]>0$ if $w\neq0$, which yields that
 $f_1(t)=-g(w_B(t))e_{1B}>0$  if $w\neq0$.
This and the fact that $2\rho_1^++\sigma>0$, we conclude that $\xi_1>0$.

We complete the proof of this lemma.

\begin{lem}\label{lemm18}
Let $q$ be JL-critical, then $\xi_1>0$.
\end{lem}
\emph{Proof.} Since (\ref{e50}) can be written as $(t^2y'_1)'=te^{-\rho_1t}f_1$, integrating on $(1,t)$, we obtain
$$
ty'_1(t)=t^{-1}y'_1(0)+t^{-1}\int_1^t e^{-\rho_1s}f_1(s)ds,~~t>1,
$$
which gives that $\lim_{t\rightarrow\infty}ty'_1(t)=0$.
As a consequence, since $\lim_{t\rightarrow +\infty}y_1(t)=\xi_1$, integrating  (\ref{e50}) on $(t,\infty)$, we obtain
$$
-(ty'_1(t)+y_1(t))+\xi_1=\int_t^\infty e^{-\rho_1s}f_1(s)ds.
$$
Since $z_1(t)=te^{\rho_1t}y_1(t)$, we have $ty'_1(t)+y_1(t)=(-\rho_1z_1(t)+z'_1(t))e^{-\rho_1t}$.
Taking $t\rightarrow0$, we have
\begin{equation}\label{e54}
\xi_1=(-\rho_1z_1(0)+z'_1(0))+\int_0^\infty e^{-\rho_1s}f_1(s)ds.
\end{equation}
Similarly, integrating  (\ref{e50}) on $(-\infty,t)$, we obtain
$$
ty'_1(t)+y_1(t)=\int_{-\infty}^t e^{-\rho_1s}f_1(s)ds,
$$
where we have used the facts that $\lim_{t\rightarrow-\infty}ty'_1(t)=0$ and $\lim_{t\rightarrow-\infty}y_1(t)=0$.
Since  $ty'_1(t)+y_1(t)=(-\rho_1z_1(t)+z'_1(t))e^{-\rho_1t}$,
taking $t\rightarrow0$, we have
\begin{equation}\label{e55}
(-\rho_1z_1(0)+z'_1(0))=\int_{-\infty}^0 e^{-\rho_1s}f_1(s)ds.
\end{equation}
Add (\ref{e54}) to (\ref{e55}), we conclude that
$$
\xi_1=\int_{-\infty}^\infty e^{-\rho_1s}f_1(s)ds,
$$
which yields $\xi_1>0$.
The proof of this lemma is complete.

For the JL-subcritical case, if $(\xi_1,\xi_2)=0$, where $(\xi_1,\xi_2)$ is given  in (\ref{e48}),  then we will establish the following  expansion with a smaller error.

\begin{lem}\label{lemm19}
Let $q$ be JL-subcritical and $(\xi_1,\xi_2)=0$, then there exists $\xi<0$ and $C,~\varepsilon>0$ such that
$$
\|e^{-\rho_2^-t} w(t,\theta)-\xi e_2\|_{2,a}\leq Ce^{-\varepsilon t}.
$$
Moreover the norm $\| \cdot\|_{2,a}$ can be replaced by the norm $\| \cdot\|_{\infty}$.
\end{lem}
\emph{Proof.}
Since $(\xi_1,\xi_2)=0$, from (\ref{e47}) and (\ref{e52}), we have $|z_1(t)|\leq Ce^{3\rho_2^-t/2}$.
Then, by the same argument as in the proof of Lemma \ref{lemm15}, we obtain
\begin{equation}\label{e56}
\sum_{i=3}^\infty |z_i(t)|\leq C[e^{\rho_3^-t}+e^{3\rho_2^-t/2}],~~~ |z_2(t)|\leq Ce^{\rho_2^-t}.
\end{equation}
Now we set $y_2(t)=e^{-\rho_2^-t}z_2(t)$. Then $y_2(t)$ solves
\begin{equation}\label{e57}
y''_2+(2\rho_2^-+\sigma)y'_1=e^{-\rho_2^-t}f_2.
\end{equation}
By (\ref{e56}) we know that $y_2(t)$ is bounded, so
there exists a sequence $\{t_i\}_{i\in N}$ such that
$
\lim_{i\rightarrow\infty}y'_2(t_i)=0.
$
Then we have
$$
-e^{(2\rho_2^-+\sigma)t}y'_2(t)=\lim_{i\rightarrow\infty}e^{(2\rho_2^-+\sigma)t_i}y'_2(t_i)-e^{(2\rho_2^-+\sigma)t}y'_2(t)
=\int_t^\infty e^{(\rho_2^-+\sigma)s}f_2(s)ds.
$$
Hence from (\ref{e56}), we deduce $|y'_2(t)|\leq Ce^{\rho_2^-t/2}$.
As a consequence, there exists $\xi\in \mathbb{R}$ such that $
\lim_{t\rightarrow\infty}y_2(t_i)=\xi$ and
 \begin{eqnarray}
\nonumber
\|e^{-\rho_2^-t} w(t,\cdot)-\xi e_2\|_{2,a}&\leq& \|e^{-\rho_2^-t} w(t,\cdot)-e^{-\rho_2^-t}z_2(t) e_2\|_{2,a}+
\|(e^{-\rho_2^-t}z_2(t)-\xi) e_2\|_{2,a}
\\
\nonumber &
\leq&\|e^{-\rho_2^-t}[ w(t,\cdot)-z_2(t) e_2]\|_{2,a}+\|(y_2(t)-\xi) e_2\|_{2,a}
\\
\nonumber &
\leq& Ce^{-\varepsilon t},
\end{eqnarray}
where in the last inequality we have used  (\ref{e56}), the facts that $|z_1(t)|\leq Ce^{3\rho_2^-t/2}$ and exponent decay of $y'_2(t)$.
Integrating (\ref{e57} on $(-\infty,\infty$), we obtain
$$
(2\rho_2^-+\sigma)\xi=\int_{-\infty}^\infty e^{-\rho_2^-s}f_2(s)ds.
$$
since $2\rho_2^-+\sigma<0$ and $f_2(s)=-g(w_B(t))e_{2B}>0$, we deduce that $\xi<0$.
Applying the same argument as given in Lemma  \ref{lemm16}, we obtain the $L^\infty$-estimate.

We complete the proof of this lemma.

\emph{Proof of Proposition \ref{prop1}}
By Lemmas \ref{lemm16}-\ref{lemm19}, we immediately obtain Proposition \ref{prop1}.

\vskip0.1in

{\bf Acknowledgment.} The author is partially supported by  NSFC,  No   11101134 and NSFC, No 11371128.

\end{document}